\documentclass[12pt]{amsart}
\usepackage{amscd,amssymb}
\usepackage{tikz}
\usepackage[plainpages,backref,urlcolor=blue]{hyperref}
\usepackage[all]{xy}

\theoremstyle{plain}
\newtheorem{thm}[subsection]{Theorem}
\newtheorem{lem}[subsection]{Lemma}
\newtheorem{prop}[subsection]{Proposition}
\newtheorem{cor}[subsection]{Corollary}

\theoremstyle{definition}
\newtheorem{rk}[subsection]{Remark}
\newtheorem{definition}[subsection]{Definition}
\newtheorem{ex}[subsection]{Example}
\newtheorem{conj}[subsection]{Conjecture}

\numberwithin{equation}{section}
\newcommand{\OO}{{\mathcal O}}

\newcommand{\wI}{\widehat{I}}

\newcommand{\CC}{{\mathcal C}}

\newcommand{\Z}{\mathbb{Z}}

\newcommand{\C}{\mathbb{C}}
\newcommand{\PP}{\mathbb{P}}

\newcommand{\N}{\mathbb{N}}

\begin{document}

\title [Lefschetz properties of Jacobian algebras and Jacobian modules]
{Lefschetz properties of Jacobian algebras and Jacobian modules}

\author[Alexandru Dimca]{Alexandru Dimca$^{1}$}
\address{Universit\'e C\^ ote d'Azur, CNRS, LJAD, France and Simion Stoilow Institute of Mathematics,
P.O. Box 1-764, RO-014700 Bucharest, Romania}
\email{dimca@unice.fr}

\author[Giovanna Ilardi]{Giovanna Ilardi}
\address{Dipartimento Matematica Ed Applicazioni ``R. Caccioppoli''
Universit\`{a} Degli Studi Di Napoli ``Federico II'' Via Cintia -
Complesso Universitario Di Monte S. Angelo 80126 - Napoli - Italia}
\email{giovanna.ilardi@unina.it}


\thanks{\vskip0\baselineskip
\vskip-\baselineskip
\noindent $^1$This work has been partially supported by the Romanian Ministry of Research and Innovation, CNCS - UEFISCDI, grant PN-III-P4-ID-PCE-2020-0029, within PNCDI III}



\begin{abstract} Let $V:f=0$ be a  hypersurface of degree $d \geq 3$ in the complex projective space $\PP^n$, $n \geq 3$, having only isolated singularities. Let $M(f)$ be the associated Jacobian algebra and $H: \ell=0$ be a hyperplane in $\PP^n$ avoiding the singularities of $V$, but such that $V \cap H$ is singular. We related the Lefschetz type properties of the  linear maps $\ell: M(f)_k \to M(f)_{k+1}$ induced by the
multiplication by linear form $\ell$  to the singularities of the hyperplane section $V \cap H$. Similar results are obtained for the Jacobian module $N(f)$.

\end{abstract}
 
\maketitle


\section{Introduction} 

Let $S=\C[x_0,\ldots,x_n]$ be the graded polynomial ring in $n+1$ variables with complex coefficients, with $n \geq 2$. We denote by $S_k$ the degree $k$ homogeneous component of $S$.
Let $f \in S_d$ be a homogeneous polynomial such that the hypersurface $V:f=0$ in the projective space $\PP^n$ is smooth. Then the Jacobian algebra $M(f)=S/J(f)$ is a standard graded Artinian Gorenstein algebra with socle degree $T=(n+1)(d-2)$, where we set
$$f_j= \frac{\partial f}{\partial x_j}$$
for $j=0, \ldots,n$ and $J(f)=(f_0, \ldots,f_n) \subset S$ is the Jacobian ideal of $f$. This algebra is also known as the Milnor algebra of $f$. In fact, the hypersurface $V:f=0$  is smooth if and only if $M(f)$ is an Artinian algebra. 

\begin{definition}
Let $M=\displaystyle{\oplus_{i=0}^T}M_i$ be an Artinian graded $\C-$algebra with $M_T\neq 0$. 
\begin{enumerate}
\item The algebra $M$ is said to have the Weak Lefschetz Property in degree $i$, for $i < T$, for short $WLP_i$, if there exists an element $L\in M_1$ such that the multiplication map $ L: M_i\to M_{i+1}$ is of maximal rank. 
We say that the algebra $M$ has $WLP$ if it has $WLP_i$  for all $0\leq i \leq T-1$.

\item  We say that $M$ has the Strong Lefschetz Property in degree $k<T/2$, for short $SLP_k$, if there is $L \in M_1$ such that the linear map $ L^{T-2k}: M_k \to M_{T-k}$ is an isomorphism. We say that the algebra $M$ has $SLP$  if it has $SLP_k$ in degree $k$  for all $k <T/2$. 
\end{enumerate}
\end{definition}

\begin{ex}\label{ex:monomial} This example is due to Stanley  \cite{St} and  Watanabe  \cite{Wa3} and it is considered to be the starting point of the research area of Lefschetz properties for graded algebras. 
Consider the standard graded Artinian Gorenstein algebra
 $$M= \frac{\C[x_0,\ldots,x_n]}{(x_0^{a_0},\ldots,x_n^{a_n})}$$
  with integers $a_i>0$, for all $i=0,\ldots,n$. 
 Then $M$ has $SLP$.
\end{ex}
In particular, for a Fermat type polynomial
$$f_F=x_0^d+ \ldots + x_n^d,$$
the Jacobian algebra $M(f_F)$ has $SLP$. By obvious semi-continuity properties, it follows that the Jacobian algebra $M(f)$ has $SLP$ for a {\it generic } polynomial $f \in S_d$. However, it seems that we have no control on the meaning of the word 'generic' in this claim.
In fact, one has the following.

\begin{conj}\label{conj1}
The Jacobian algebra $M(f)$ has $SLP$ for any polynomial $f \in S_d$ such that the associated hypersurface $V:f=0$ is smooth.
\end{conj}

\begin{rk}\label{rk1}
We list here the known  results related to this conjecture.

\begin{enumerate}
\item
When $n=2$,
for any smooth curve $V:f=0$, the associated Jacobian algebra $M(f)$ has  $WLP$, as follows from the more general results in \cite{HMNW}. 
Moreover, when $d=2d'$ is even, then the multiplication by the square of a generic linear form $\ell \in S_1$ induces an isomorphism
 $$\ell^2: M(f)_{3d'-4} \to M(f)_{3d'-2}.$$
In particular, when $d=4$, the Jacobian algebra $M(f)$ has the  $SLP$, see \cite{DGI}.

\item When $n=3$, if $V:f=0$ is any smooth cubic surface in $\PP^3$, then $M(f)$ has the $SLP$, see \cite{DGI}. In fact, this result  is a consequence of the classical Hesse-Gordan-Noether’s
Theorem, see for instance \cite{G,R}.
Moreover, for any smooth surface $V(f)$, the $SLP$ holds in degree 1, see for instance \cite[Theorem B]{BFP}. For any smooth surface  in $\PP^3$ of degree $d \in \{4,5,6\}$, the $WLP$ holds in all degrees, see \cite[Corollary 7.2]{B+}. Moreover, for degrees $d >6$, the $WLP$ holds in all degrees
up to $(3d-2)/2$, see \cite[Corollary 7.3]{B+}.

\item For $n=4$,  if $V:f=0$ is any smooth cubic 3-fold in $\PP^4$, then $M(f)$ has the $SLP$, \cite[Theorem C]{BFP}. 

\item For arbitrary dimension $n$ and degree $d$, one knows that  $M(f)$ has the $WLP$ in degree $\leq d-2$, see \cite{I1}. The $SLP_0$ also holds for obvious reasons. Indeed, if $\ell^T \in J(f)$ for all $\ell \in S_1$, then the hessian polynomial $hess(f) \notin J(f)$ would have no Waring decomposition.

\item In \cite{A} WLP is proven for some complete intersection algebras
presented in degree 2, improving some previous bounds. Finally, in \cite{B}, SLP1 for the Jacobian Algebra of a smooth cubic fourfold
is proved.

\end{enumerate}

\end{rk}

\begin{rk}\label{rk2}
Using the duality of the Jacobian algebra $M(f)$, in order to prove the $WLP$ for
$M(f)$, it is enough to prove the $WLP_i$ for $i <T/2$.
\end{rk}

In fact, one can define $WLP$ for more general graded objects as follows.

\begin{definition}
\label{defWLPE}
Let $M=\displaystyle{\oplus_{i\geq 0}}M_i$ be a graded $\C-$algebra (resp. a graded $S$-module) with  $\dim M_i <  + \infty$  for any $i \geq 0$.
 The algebra (resp. module) $M$ is said to have the Weak Lefschetz Property in degree $i$,  for short $WLP_i$, if there exists an element $L\in M_1$ (resp. $L \in S_1$) such that the multiplication map $ L: M_i\to M_{i+1}$ is of maximal rank. 
We say that $M$ has $WLP$ if it has $WLP_i$  for all $i \geq 0$.
\end{definition}
When the hypersurface $V:f=0$ in $\PP^n$ is singular, then one can consider  the Jacobian algebra $M(f)=S/J(f)$ and the Jacobian module $N(f) \subset M(f)$ given by $N(f)=I(f)/J(f)$, where $I(f)$ is the saturation of $J(f)$ with respect to the maximal ideal ${\bf m} =(x_0, \ldots,x_n)$, see \cite{Se}. Recall that, for any homogeneous ideal $I$ in $S$, we define its saturation $\wI$ with respect to $\bf m$ as the set of all elements $s \in S$ such that for any $i=0,...,n$ there is a positive integer $m_i$ such that 
$x_i^{m_i}s \in I$.
One has $N(f)=M(f)$ if and only if $V:f=0$ is smooth, so $N(f)$, which is always an Artinian module, can be thought of as a replacement of $M(f)$ when $V $ is singular.

\begin{ex}\label{exCurves} 
Let $C:f=0$ be a reduced degree $d$ curve in $\PP^2$. Then the corresponding Jacobian algebra $M(f)$ satisfies $WLP_i$ for any
$$i <\frac{3(d-2)}{2}$$
and the corresponding mappings $M(f)_i \to M(f)_{i+1}$ are injective.
Moreover, the Jacobian algebra $M(f)$ satisfies $WLP_i$ for any
$$i \geq i_0=  \Big{\lfloor } \frac{3(d-2)}{2} \Big {\rfloor}$$
and the corresponding mappings $M(f)_i \to M(f)_{i+1}$ are surjective if and only if 
\begin{equation}
\label{eqCT}
ct(C) \geq 3(d-2)-i_0,
\end{equation}
 see \cite[Corollary 4.4]{DP}. Here
$$ct(C)=\max\{k \in \Z \ : \ \dim M(f)_j= \dim M(f_F)_j \text{ for all } j \leq k\},$$
where $f_F$ denotes the Fermat degree $d$ polynomial as in Example
\ref{ex:monomial}. As explained in \cite[Remark 4.5]{DP}, many curves satisfy the condition \ref{eqCT}, in particular all nodal curves. Moreover,
 it has been shown in \cite[Corollary 4.3]{DP}
that the Jacobian module $N(f)$ satisfies $WLP$ for any reduced plane curve.
\end{ex}
\begin{ex}\label{exHyp} 
Let $V:f=0$ be a  degree $d$ hypersurface in $\PP^n$ having only isolated singularities. Then the corresponding Jacobian algebra $M(f)$ satisfies $WLP_i$ for any
$$i >n(d-2)$$
and the corresponding mappings $M(f)_i \to M(f)_{i+1}$ are surjective.
For $i>(n+1)(d-2)$ these mappings are isomorphisms, see \cite[Corollary 8]{CD}.
\end{ex}

In this note, we fix a  hypersurface $V:f=0$  in $\PP^n$ having only isolated singularities and a hyperplane $H: \ell=0$ in $\PP^n$ avoiding the singularities of $V$, with $n , d \geq 3$.   The hyperplane section $V(f, \ell)=V(f) \cap H$ has then at most isolated singularities, as follows from the following.
\begin{lem}
\label{lem1}
Let $V$ be a hypersurface in $\PP^n$ of degree $d$, with $n,d \geq 2$. Assume that $H$ is a hyperplane such that $H \cap V_s=\emptyset$, where $V_s$ denotes the singular set of $V$. Then $V$ and $V \cap H$ have at most isolated singularities. Moreover, when $V$ has at most isolated singularities, then there are hyperplanes $H$ such that $H \cap V_s=\emptyset$ and $V \cap H$ is singular, except when $V$ is a cone.
\end{lem}
As usual, $V$ is a cone if  one can choose the coordinates $x$ such that $f_0=0$, hence $f$ does not depend on the variable $x_0$.
If $V:f=0$ is a cone with an isolated singularity, then up to a linear change of coordinates, we have
that $f$ is a polynomial $g\in \C[x_1,...,x_n]$  and the corresponding hypersurface $g=0$ is smooth in $H:x_0=0$.
Then obviously the multiplication by $x_0$ is injective in any range, since
$$M(f)=M(g) \otimes \C[x_0].$$

\begin{rk}\label{rk2.5}
 Note that any hypersurface $W \subset H$ with only isolated singularities may occur as a section 
$W=V \cap H$ for a certain smooth hypersurface $V$, see \cite[Proposition (11.6)]{RCS}.
\end{rk}

In this setting, and under the assumption that $V(f, \ell)$ is singular, we investigate the injectivity of the multiplication maps $\ell :M(f)_k \to M(f)_{k+1}$. 
Our  main result is stated in terms of some numerical invariants of the hyperplane section $V(f, \ell)$, which we define now. Clearly $V(f, \ell)$ is a hypersurface in $H=\PP^{n-1}$. If we choose a system of coordinates $y=(y_1, \dots,y_n)$ on $H$, then  $V(f, \ell)$ given by an equation $g=0$, hence it has a Jacobian ideal
$J(g)$ in the polynomial ring $R=\C[y_1,...,y_n]$. Let $I(g)$ be the saturation of the Jacobian ideal
$J(g)$ with respect to the maximal ideal $(y_1,...,y_n)$ and let
$s(g)$ be the initial degree of the graded ideal $I(g)$, namely
\begin{equation} \label{e1} 
s(g)=\min \{j \in \N \ : \ I(g)_j \ne 0 \} \leq d-1,
\end{equation}
where $d=\deg f =\deg g.$
Let $g_j$ denote the partial derivative of $g$ with respect to $y_j$ and consider the graded $R$-module $Syz(g)$ of first order syzygies of $g_1, \ldots,g_n$, namely
\begin{equation} \label{e2} 
Syz(g)= \{a=(a_1,\ldots,a_n) \in R^n \ : \ a_1g_1+ \ldots +a_ng_n=0 \}.
\end{equation}
Let $r(g)$ be the initial degree of the graded module $Syz(g)$, namely
\begin{equation} \label{e3} 
r(g)=\min \{j \in \N \ : \ Syz(g)_j \ne 0 \} \leq d-1.
\end{equation}
It is clear that both invariants $r(g)$ and $s(g)$ do not depend on the choice of the linear coordinates $y$ on $H=\PP^{n-1}$. Note that $V(f, \ell)$ singular implies $s(g)>0$. On the other hand, $r(g)=0$ if and only if $V(f, \ell)$  is a cone. With this notation, our  main result is the following improvement of the second author result recalled in Remark \ref{rk1} (4).

\begin{thm}\label{thm1}
Let $V:f=0$ be a hypersurface in $\PP^n$ of degree $d$, with $n,d \geq 2$. Assume that $H: \ell=0$ is a hyperplane such that $H \cap V_s=\emptyset$, where $V_s$ denotes the singular set of $V$, and that $V(f, \ell)=V \cap H$ is singular. 
Then the multiplication map $\ell :M(f)_k \to M(f)_{k+1}$ is injective for any
$$k \leq  \min\{d-3 + r(g), d-3+s(g)\}.$$
\end{thm}
  Theorem \ref{thm1} can give the $WLP$  for the Jacobian algebra $M(f)$ (in some degrees or in all degrees) for some  generic classes of hypersurfaces, {\it with a precise meaning of the word 'generic'}, as well as for many new {\it non generic classes} of smooth hypersurfaces.

In terms of the Jaconian module $N(f)$, this result implies the following.
\begin{cor}\label{corN}
Let $V:f=0$ be a hypersurface in $\PP^n$ of degree $d$, with $n,d \geq 3$. Assume that $H: \ell=0$ is a hyperplane such that $H \cap V_s=\emptyset$, where $V_s$ denotes the singular set of $V$, and that $V(f, \ell)=V \cap H$ is singular. 
Then the multiplication map $\ell :N(f)_k \to N(f)_{k+1}$ is injective for any
$$k \leq  k_0=\min\{d-3 + r(g), d-3+s(g)\}.$$
Moreover, the multiplication map $\ell :N(f)_k \to N(f)_{k+1}$ is surjective for any $$k \geq (n+1)(d-2)-k_0-1.$$
\end{cor}

In section 2 we prove Lemma \ref{lem1} , Theorem \ref{thm1}, and Corollary \ref{corN}.
In section 3 we apply Theorem \ref{thm1} to surfaces $S$ in $\PP^3$ with isolated singularities. The surfaces of degree 4 are considered in Example \ref{ex1} and Example \ref{ex2}.
The surfaces having a section $C=S \cap H$ which is a nodal curve and such that all irreducible components of $C$ are rational (resp. $C$ is a free or a nearly free curve) are considered in Proposition \ref{prop1} (resp. Proposition \ref{prop2}). When the surface $S$ is smooth, all the results in this section except Proposition \ref{prop1} are obtained in stronger forms in \cite{B+}.

 In section 4, we consider higher dimensional hypersurfaces $V$, having a nodal hyperplane section $Y=V \cap H$ with many singularities, e.g. $Y$ is a Kummer surface in Proposition \ref{prop3}, respectively a Chebyshev hypersurface in Proposition \ref{prop5}.
Finally, an application to $WLP$ of a recent result of the authors, see Theorem \ref{thmDuals} below quoted from \cite{Duals},   is given in Corollary \ref{cor20} and  Corollary \ref{cor21}.
We state all the results in this section for a smooth hypersurface $V$,
but all of them, except Theorem \ref{thmDuals}, have obvious versions for a hypersurface $V$ having only isolated singularities.

\medskip

We would like to thank Arnaud Beauville for very useful discussions related to this paper. We  also  thank the referee for the very careful reading of our manuscript and for his suggestions to improve the presentation.

\section{The proofs of Lemma \ref{lem1}, Theorem \ref{thm1} and  Corollary \ref{corN}} 

\subsection{The proof of Lemma \ref{lem1} }
Since $H \cap V_s= \emptyset$, it is clear that $\dim V_s \leq 0$.
Assume that $V:f=0$ and $H:x_0=0$. The singular set of $V \cap H$ is then given by
$$Z: f_1(0,x_1,\ldots,x_n)=f_2(0,x_1,\ldots,x_n)=\ldots =f_n(0,x_1,\ldots,x_n)=0.$$
If $\dim Z >0$, then any point in the non-empty intersection of $Z$ with the hypersurface $f_0(0,x_1,\ldots,x_n)=0$ would give rise to a point in
$V_s \cap H$. This contradiction shows that $\dim Z \leq 0$.
 
 The second claim in Lemma \ref{lem1} is more subtle. Assume from now on that $V$ is not a cone.
For any singular point $p \in V_s$, denote by $H_p$ the hyperplane in the dual space $(\PP^n)^{\vee}$ corresponding to all the hyperplanes in $\PP^n$ passing through $p$. Note that the dual hypersurface $V^{\vee}$ is irreducible, since $V$ is so since $n \geq3$. It follows that if $V^{\vee}$ is contained in the union of all hyperplanes $H_p$, for $p \in V_s$, there is a point $q \in V_s$ such that $V^{\vee} \subset H_q$. We can assume that $q=(1:0: \ldots :0)$, and then this condition means
$$f_0(x)=0 \text{ for all } x \in V \setminus V_s.$$
But this implies that $f_0$ is divisible by $f$, which is possible only if $f_0=0$, that is when $V$ is a cone. This contradiction shows that there is a hyperplane $H : \ell =0$ such that $H \cap V_s=0$ and
$H$ is tangent to $V$. In fact $H$ corresponds to a point in $V^{\vee} $ not in any $H_p$, for $p \in V_s$.
This implies that the section $V\cap H$ is singular. 

\subsection{The proof of Theorem \ref{thm1} }

Without loss of generality, we may take $\ell=x_0$ and $y_j=x_j$ for all $j=1, \ldots,n$. Then one may write
$$f(x_0,y)=g(y)+x_0h(y)+x_0^2p_2(y)+ \ldots + x_0^dp_d(y),$$
where $p_j \in R_{d-j}$. Assume that $\ell :M(f)_k \to M(f)_{k+1}$ is not injective.
Then there is a homogeneous polynomial $q  \in S_k$   such that
$$q \notin J(f)_k \text{ and } x_0q \in J(f)_{k+1}.$$
It follows that one has a relation
$$x_0q=b_0f_0+ \ldots + b_nf_n,$$
with homogeneous polynomials $b_j \in S_{k+2-d}$, not all divisible by $x_0$. If we set $x_0=0$, this yields
\begin{equation} \label{e4} 
c_0h+c_1g_1+ \ldots + c_ng_n=0,
\end{equation}
where not all the polynomials $c_j(y)=b_j(0,y) \in R_{k+2-d}$ are zero.
There are two cases to discuss.

\medskip

\noindent {\bf Case 1.} $c_0 \ne 0$. Then note that the singular set of the hypersurface $V$  on the hyperplane $x_0=0$ is given by the solutions of the system
$$h=g_1= \ldots = g_n=0.$$
Since $V_s \cap H= \emptyset$,  and $V(f, \ell)$ is assumed to be singular, it follows that $h$ does not vanish on the singular set of the section $V(f, \ell)$, which is given by the solutions of the system
$$g_1= \ldots = g_n=0.$$
In other words, if $q$ is a singular point of $V(f, \ell)$, and $u \in R_{d-1}
$ is a homogeneous polynomial not vanishing at $q$, the regular function germ of $h/u$ at $q$ is an invertible element in the corresponding local ring $\OO_{H,q}$.
This observation and the equation \eqref{e4} imply that $c_0 \in I(g)$, see the discussion at the beginning of Section 2 in \cite{DBull}. Indeed,  equation \eqref{e4} implies that the germ of rational function
associated to $c_0$ at each singular point $q$ of $V(f, \ell)$ belongs to the stalk at $p$ of the sheaf associated to the Jacobian ideal at $g$, which is an ideal in the local ring $\OO_{H,q}$.
This property implies that $c_0 \in I(g)$.
Hence 
$ k+2-d =\deg c_0 \geq s(g)$ which gives $k \geq d-2+s(g)$.
\medskip

\noindent {\bf Case 2.} $c_0 = 0$.
Then the equation \eqref{e4} becomes a non-zero homogeneous element  $\rho \in Syz(g)$. It follows that $ k+2-d =\deg \rho \geq r(g)$. In other words, we have $k \geq d-2+r(g)$.

\medskip

This ends the proof of Theorem \ref{thm1}.

\subsection{The proof of Corollary \ref{corN}}

The first claim follows from Theorem \ref{thm1}, since $N(f)$ is a submodule in $M(f)$. The second claim follows from the duality property of the Jacobian module $N(f)$, namely
$$N(f)_j^{\vee} \cong N(f)_{T-j},$$
for any integer $j=0, \ldots, T$, where $T=(n+1)(d-2)$, see either \cite[Theorem 3.4]{Se} or \cite[Theorem 4.7]{vSt}.

\section{Surfaces of degree $d \geq 4$ in $\PP^3$} 
In this case the plane section $V(f, \ell)$ is a plane curve of degree $d$, with isolated singularities. We discuss several possibilities.

\begin{ex}\label{ex1}
For smooth surfaces of degree 4 $WLP$ holds in all degrees, see \cite[Corollary 7.2]{B+}.
Let $V:f=0$ be a  surface of degree 4 in $\PP ^3$ having isolated singularities and containing no line, and let $H: \ell=0$
be a plane such that  $V(f)_s \cap H=\emptyset$  and  the plane curve $V(f, \ell)=V \cap H$ has 
 at least 3 singular points with total Tjurina number $\tau(V(f,\ell))$ at most 5.
Then clearly $s(g) \geq 2$, since any element in $I(g)$ vanishes at the singular points  and
$V(f, \ell)$ has no line components. Using the lower bound on $\tau(V(f,\ell))$ given in \cite{dPW}, we see that $r(g) \geq 2$. It follows from Theorem \ref{thm1}  and Corollary \ref{corN} that
the multiplication map $\ell :M(f)_i \to M(f)_{i+1}$ is injective for $i \leq 3$ and that the Jacobian module $N(f)$ has $WLP$.  Other quartic surfaces are discussed in Example \ref{ex2}. 
\end{ex}

\begin{prop}\label{prop1}
If the surface $V:f=0$ with isolated singularities admits a plane section $V(f,\ell)=V \cap H$, where $H: \ell=0$, such that $V_s \cap H=\emptyset$  and $V(f,\ell)$ is a nodal curve, all irreducible components of $V(f,\ell)$ being rational curves, then the corresponding Jacobian algebra $M(f)$ has the $WLP_k$ for
$k<2d-4$ and the Jacobian module $N(f)$ has $WLP$.
\end{prop}
\proof
 Using \cite[Theorem 4.1]{DStEdin}, it follows that $r(g)\geq d-2$. On the other hand, we have $s(g) \geq d-2$, using 
\cite[Theorem 3.2]{CDI}. 
It follows from Theorem \ref{thm1} that
the multiplication map $\ell :M(f)_k \to M(f)_{k+1}$ is injective for
$k <2d-4=T/2$. We conclude  using Corollary \ref{corN}.
\endproof

\begin{rk}\label{rk3}
Any generic quartic surface in $\PP^3$ admits a section which is a rational nodal curve, see \cite[Theorem 1.2]{Chen1}. 
On the other hand, for $d\geq 5$, a generic degree $d$ surface in $\PP^3$ does not admit a section which is a rational nodal curve. More precisely, in this case any irreducible component $C$ of a section $V(f, \ell)$ of a generic surface has geometric genus satisfying the inequality
$$g(C) > \frac{d(d-3)}{2}-3,$$
see  \cite[Theorem 1]{Xu}. In the same paper, Xu shows that the list of singularities on the section $V(f, \ell)$ is one of the following, see  \cite[Proposition 3]{Xu}. 
\begin{enumerate}
\item $A_1$, $2A_1$, $3A_1$;

\item $A_2$, $A_1A_2$

\item $A_3$.
\end{enumerate}

In other words, only list of singularities with total Tjurina number $\tau(V(f, \ell)) \leq 3$ may occur for a generic surface of degree $d \geq 5$. For a similar result, see  also \cite{Bruce}.

\end{rk}

For a curve $C:g=0$, we set $n(f)_k=\dim N(g)_k$ for any integer $k$ and also $\nu(C)=\max _j \{n(g)_j\}$.
Recall that the curve $C$ is free if and only if $N(g)=0$ and hence $\nu(C)=0$, and $C$ is nearly free if and only if $\nu(C)=1$, see for instance \cite{CDI, DStRIMS}. For a free (resp. nearly free) curve, the exponents are
$d_1=r(g)$, the initial degree of the graded $R$-module $Syz(g)$
and $d_2=d-1-d_1$ (resp. $d_2=d-d_1$).

\begin{prop}\label{prop2}
If the surface $V:f=0$ with isolated singularities admits a plane section $V(f,\ell)=V \cap H$, where $H: \ell=0$, such that $V_s \cap H=\emptyset$  and $V(f,\ell)$ is a free (resp. nearly free) curve with exponents $(d_1,d_2)$, then the multiplication map $\ell :M(f)_k \to M(f)_{k+1}$ is injective for any
$$k \leq d+d_1-3.$$
\end{prop}
In the case of smooth surfaces $V$, this result is weaker than \cite[Corollary 7.2]{B+}.
\proof
For a free curve $V(f, \ell):g=0$ one has $J(g)=I(g)$ and hence
$s(g)=d-1$. By definition one has $r(g)=d_1$ and it is known that
$2d_1 \leq d-1$ which implies $d_1 \leq d-1=s(g)$. This proves our claim for a free curve.

For a nearly free curve $V(f, \ell):g=0$ one has $s(g)=\min \{d-1, d+d_1-3\}$, see \cite[Corollary 2.17]{DStRIMS}. By definition 
$r(g)=d_1$ and it is known that
$d_1+d_2=d$. It follows that $d_1 \leq d-1$ and also
$d_1 \leq d+d_1-3$ since we suppose $d \geq 3$.
\endproof
In view of the point (4) in Remark \ref{rk1}, the above result is useful only when $d_1 \geq 2$.
\begin{ex}\label{ex2}
Let $V$ be a surface of degree 4 with isolated singularities in $\PP ^3$  and let $H: \ell=0$
be a plane avoiding the singularities of $V(f)$ and such that the plane curve $V(f, \ell)=V(f) \cap H$ has one of the following lists of singularities, not covered by the discussion in Example \ref{ex1}.
\begin{enumerate}
\item $3A_2$;

\item $A_2A_4$

\item $A_6$

\item $6A_1$

\item $A_1A_5$.

\end{enumerate}
In the first 3 cases the curve $V(f, \ell)$ is irreducible and rational, while in case (4) the curve $V(f, \ell)$ is a union of 4 lines in general position.
In particular, the case (4) can be treated using also Proposition \ref{prop1}. In case (5), the curve $V(f, \ell)$ is a union of two conics, meeting at 2 points, with intersection multiplicities 1, and respectively 3.
All these  curves are shown to be nearly free with exponents $(2,2)$, see
\cite[Example 2.13]{DStRIMS} and \cite[Proposition 5.5]{DIS}.
 It follows from Proposition \ref{prop2} that
the multiplication map $\ell :M(f)_k \to M(f)_{k+1}$ is injective for $k \leq 3$.
\end{ex}
\begin{rk}\label{rk4}
All the possibilities of the singularities of a plane section of a smooth quartic surface in $\PP^3$ are listed in \cite{BG}. 
We have discussed above only some of these possibilities, assuming that they occur also as sections of quartic surfaces with isolated singularities.\end{rk}

\section{Hypersurfaces  in $\PP^n$, $n \geq 4$, having sections with many nodes} 

As already said in Introduction, all the results in this section are stated for a smooth hypersurface $V$,
but all of them, except Theorem \ref{thmDuals}, have obvious versions for a hypersurface $V$ having only isolated singularities.
In this section the hyperplane section $V(f, \ell) \in H=\PP^{n-1}$ is assumed to have only nodes, namely ordinary double points, also known as $A_1$-singularities. In this case, it is known that
\begin{equation} \label{e10} 
r(g)=d-1,
\end{equation}
when $d \geq 3$, see \cite{DFerr2017}. To get lower bounds on $s(g)$ we impose a large number of singularities on the hyperplane section $V(f, \ell)$.
We discuss several cases.

\begin{prop}\label{prop3}
If the smooth 3-fold $V:f=0$ in $\PP^4$ of degree $4$ admits a hyperplane section $V(f, \ell):g=0$ that is a nodal surface, with at least 10 nodes, not all on a quadric, then the corresponding Jacobian algebra $M(f)$ has the $WLP$. In particular, this occurs when the
hyperplane section $K=V(f, \ell):g=0$  is a Kummer surface with 16 nodes.
\end{prop}
\proof
Using the equality \eqref{e10}, it follows that $r(g) =3$. In the first claim, we assume that $s(g) \geq 3$. This claim follows now from Theorem \ref{thm1} and Remark \ref{rk2}, since $T=10$ in this case.

 We show now that $s(g) \geq 3$ in the case of a section being a Kummer surface.
 Let $X$ be the minimal resolution of $K$. Denote by $H'$ the pull-back on $X$ of a plane section of $K$ and let $E_1, \ldots, E_{16}$ be the exceptional divisors. Assume there is a quadric $Q$ in the hyperplane $H=\PP^3$, passing through all the 16 nodes of $K$. Then the pull-back on $X$ of this quadric $Q$ gives an effective divisor $D$ in the linear system 
$$|2H'-\sum_{i=1}^{16}E_i|.$$
On the other hand, there are 16 planes in $H$, classically called tropes, each of them tangent to $K$ along a conic passing through 6 of the singularities of $K$, see \cite[Chapter 1]{Hu}. Hence the proper transform $C_j$ of one of these conics satisfies
$$2C_j=H'- \sum_{i \in I_j}E_i,$$
with $|I_j|=6$. It follows that
$$C_j \cdot D=\frac{1}{2}(H'- \sum_{i \in I_j}E_i)\cdot (2H'-\sum_{i=1}^{16}E_i)=-2.$$
Since $C_j$ is irreducible, this implies that $D$ contains the curve $C_j$. Note that $D \cdot H'=8$, $C_j \cdot H'=2$ for all $j=1, \ldots ,16$, hand hence $D$ cannot contain all the conics $C_j$. This contradiction shows that $I(g)_2=0$, and hence $s(g) \geq 3$.
\endproof

There is a family of nodal hypersurfaces $X:g=0$, in any dimension $n$ and degree $d$, for which the subtle invariant $s(g)$ is known. They are the Chebyshev hypersurfaces $\CC(n,d) \subset \PP^n$, defined as follows. For more details see \cite{DStAlex}. Consider the $d$-th Chebyshev polynomial 
$$T_d(x)=\cos (d \arccos (x)).$$
Then the affine part of Chebyshev hypersurface $\CC(n,d)$, that is before homogenization with respect to $x_0$, is defined  by the affine equation
$$g(n,d)=T_d(x_1)+ \dots + T_d(x_n)=0$$
when $n$ is even, and by 
$$g(n,d)=T_d(x_1)+ \dots + T_d(x_n)+1=0$$
when $n$ is odd. The equality
\begin{equation} \label{e11} 
s(g)=d-2,
\end{equation}
follows from  \cite[Proposition 3.1]{DStAlex}.

\begin{prop}\label{prop5}
If the smooth hypersurface $V:f=0$  in $\PP^n$, $n \geq 4$, of degree $d \geq 3$ admits a hyperplane section $V(f, \ell):g=0$ which  is a Chebyshev hypersurface, then the multiplication map $\ell :M(f)_k \to M(f)_{k+1}$ is injective for any
$$k \leq 2d-5.$$
\end{prop}
\proof
It is enough to use Theorem \ref{thm1} and the formulas \eqref{e10} and 
\eqref{e11}.
\endproof
We recall the following result, see \cite{Duals}.

\begin{thm}\label{thmDuals}
Generic hypersurfaces  of degree $d$ in $\PP^n$, with $d,n \geq 3$, have hyperplane sections with $n$ nodes in general position. 
\end{thm}
Using this result, we can improve by one, in the case of generic hypersurfaces, the result of the second author, mentioned above in Remark \ref{rk1}, (4).

\begin{cor}\label{cor20}
If the smooth hypersurface $V:f=0$ in $\PP^n$, $n \geq 3$, of degree $d \geq 3$ admits a  nodal hyperplane section $V(f,\ell)=V \cap H$, which has exactly $n$ singularities in general position,
then the multiplication map $\ell :M(f)_k \to M(f)_{k+1}$ is injective for any
$k \leq d-1.$
In particular, this property holds for a generic hypersurface $V$ in $\PP^n$, $n \geq 3$, of degree $d \geq 3$. 
\end{cor}
\proof
It is enough to use Theorem \ref{thm1}, the formulas \eqref{e10} and 
the obvious fact that $s(g) \geq 2$ in this situation, since the nodes are in general position.

\endproof
When $n=3$, Corollary \ref{cor20} is much weaker than \cite[Corollary 7.3]{B+}.
The following is a special case of the previous result, for $n=5$.
\begin{cor}\label{cor21}
If the smooth 4-fold $V:f=0$ in $\PP^5$ of degree $3$ admits a hyperplane section $V(f, \ell)=V \cap H:g=0$ that is a nodal 3-fold, with at least 5 nodes, not all on a hyperplane in $H$, then the corresponding Jacobian algebra $M(f)$ has the $WLP$. In particular, this occurs for a generic 4-fold $V(f):f=0$ in $\PP^5$ of degree $3$ and for a smooth 4-fold $V:f=0$ in $\PP^5$ of degree $3$ having 
a hyperplane section $V(f, \ell):g=0$ which is a Segre 3-fold with 10 nodes.

 \end{cor}

\proof
Using the inequality \eqref{e10}, it follows that $r(g) \geq 2$. The assumptions in the first claim imply that $s(g)\geq 2$. This claim follows now from Theorem \ref{thm1} and Remark \ref{rk2}, since $T=6$ in this case. The fact that a  this situation occurs for a generic 4-fold $V(f):f=0$ in $\PP^5$ of degree $3$ follows from Corollary \ref{cor20}

We treat next the case of the Segre 3-fold section. The Segre 3-fold  is unique up-to a projective transformation, and can be given by the equation
$$X:g=x_1^3+x_2^3+x_3^2+x_4^3+x_5^3-(x_1+x_2+x_3+x_4+x_5)^3=0,$$
where $x_1,x_2,x_3,x_4,x_5$ are the coordinates on $H=\PP^4$.
The 10 nodes are located at $(1:1:1:-1:-1)$ and the other 9 points obtained by permuting the coordinates. Using this description of the singular set of $X$, it follows that $I(g)_1=0$, and hence $s(g) \geq 2$.
We conclude as for the first claim above.
\endproof

\end{document}